\documentclass[12pt,a4paper]{article}
\usepackage{a4wide}
\usepackage{amsfonts,amsmath}
\usepackage{latexsym}
\allowdisplaybreaks
\usepackage{times}
\title{The isoperimetric problem for $3$-polytopes with six vertices } 

\author{ K\'aroly J. B\"or\"oczky\footnote{Supported by NKFIH projects ANN 121649, K 109789,
 K 116451 and KH 129630}, \'Agnes Kov\'acs}

\newcommand{\proof}{\noindent{\it Proof: }}
\newcommand{\proofbox}{\mbox{ $\Box$}\\}
\newcommand{\R}{\mathbb{R}}

\newtheorem{lemma}{Lemma}[section]
\newtheorem{theo}[lemma]{Theorem}

\newtheorem{prop}[lemma]{Proposition}

\begin{document}

\maketitle

\section{Introduction}

By a convex body in $\R^d$, $d\geq 2$, we mean a compact convex set $K$ with non-empty interior.
In addition $K$ is a polytope if it is the convex hull of finitely many points.
Classically, the isoperimetric problem in $\R^d$ asks for the minimum of
the surface area of a convex body with given volume; or equivalently, the
minimum of the isoperimetric ratio $S(K)^d/V(K)^{d-1}$ for a convex body $K$ where $S(K)$ is the surface area and $V(K)$ is the volume. The minimal body is
naturally the suitable ball. The solution to this problem in the planar case was
already known to the ancient Greeks. In $\R^d$ with $d\geq 3$, the first
proofs were provided with the help of Steiner's symmetrization in the middle of the 19th century
by Steiner \cite{Ste38}. Important later contributors are, among others,
Schwarz, Minkowski, Blaschke, Hadwiger, de Giorgi. By their work, the optimality of the ball has
been also verified for a much wider class of sets (see Talenti \cite{Tal93}).
For extensive survey about this subject, see Florian \cite{Flo93} up to 1993, and Schneider \cite{Sch14} for more recent developments.

For polytopes, the meaningful question is to search for a polytope of minimal isoperimetric ratio
in a given combinatorial type. 
The isoperimetric problem for polygons of given number of vertices was settled by K. Weierstrass
(cf. Mathematische Werke, Vol. 7, p. 70-77). In higher dimension, using Steiner symmetrization,
Hadwiger \cite{Had57} proved that the regular simplex has the minimal surface area among simplices of given volume
in $\R^d$, $d\geq 3$, but the result itself was probably known to Steiner. Since this paper concentrates on dimension three, we quote the result only for $d=3$.

\begin{theo}[Hadwiger]
\label{isoper-4vertices}
For any tetrahedron $P$ in $\R^3$, we have
$$
\frac{S(P)^3}{V(P)^2}\geq 216\sqrt{3}
$$
with equality if and only if $P$ is a regular tetrahedron.
\end{theo}

For $n\geq d+2$, the isoperimetric problem for $d$-dimensional polytopes of $n$ facets 
($(d-1)$-dimensional faces) have been intensively investigated. 
According to the Lindel\"of theorem \cite{Lin69}, the polytope with minimal isoperimetric ratio contains a ball tounching every facet. Using this result, Fejes T\'oth \cite{LFT48} showed that the cube is optimal among $3$-polytopes of $6$ faces, and the regular dodecahedron is optimal among $3$-polytopes of $12$ faces.

We rather search for the optimal body among polytopes of given number $n\geq d+2$
of vertices.  If a $d$-dimensional polytope has $d+2$ vertices, the isoperimetric problem was solved 
by Boroczky, Boroczky Jr. \cite{BoB96}. In dimension $3$, this reads as follows.

\begin{theo}[Boroczky, Boroczky]
\label{isoper-5vertices}
If a polytope $P$ in $\R^3$ has at most $5$ vertices, then
$$
\frac{S(P)^3}{V(P)^2}\geq 243\sqrt{2}
$$
with equality if and only if $P$ is a double pyramid over regular triangle of edge length $a>0$ where the other six edges are of length $\frac{\sqrt{6}}4\,a$.
\end{theo}

Therefore the first open case is when a $3$-polytope has $6$ vertices.
It was proved in Hadwiger \cite{Had57} that the regular octahedron is optimal in its combinatorial type
but the results was probably known to Steiner.

\begin{prop}[Hadwiger]
\label{likeoctahedron}
If a triangulation of the surface of a $3$-polytope $P$ in $\R^3$ has the same combinatorial type as an octahedron, then
$$
\frac{S(P)^3}{V(P)^2}\geq 108\sqrt{3}
$$
with equality if and only if $P$ is a regular octahedron.
\end{prop}

The main result of this paper is that the regular octahedron is optimal among $3$-polytopes of $6$ vertices.

\begin{theo}
\label{isoper-6vertices}
If a polytope $P$ in $\R^3$ has at most six vertices, then
$$
\frac{S(P)^3}{V(P)^2}\geq 108\sqrt{3}
$$
with equality if and only if $P$ is a regular octahedron.
\end{theo}

It is interesting to note that it is not known whether the regular icosahedron is optimal in its combinatorial type
(see Brass, Moser, Pach \cite {BMP96}). 

Comparing the constants of Theorems~\ref{isoper-4vertices}, \ref{isoper-5vertices}
and \ref{isoper-6vertices}, we have $216\sqrt{3}=374.12$, $243\sqrt{2}=343.65$ and $108\sqrt{3}=187.06$.
Therefore allowing five vertices instead of four does not improve too much on the isoperimetric constant involved, while allowing six vertices substantially improves the  isoperimetric constant.

For the sake of completeness, we provide a proofs of Theorem~\ref{isoper-4vertices},
Theorem~\ref{isoper-5vertices} and Proposition~\ref{likeoctahedron}, as well.
Concerning the structure of the paper, we review notion and basic propertes of Steiner symmetrization, 
and verify Theorem~\ref{isoper-4vertices} and Proposition~\ref{likeoctahedron} in Section~\ref{secSteinersym}, 
and prove Theorem~\ref{isoper-5vertices} and Theorem~\ref{isoper-6vertices} in Section~\ref{secnoocta}.

\section{Steiner symmetrization}
\label{secSteinersym}

For points $x_1,\ldots,x_n\in\R^3$, we write $[x_1,\ldots,x_n]$ to denote their convex hull. For any 
plane $L$ containing the origin in $\R^3$ and for any $X\subset \R^3$, we write $X|L$ to denote the orthogonal projection of $X$ into $L$. For $v\in\R^3\backslash\{o\}$, we write $v^\bot$ to denote the linear $2$-plane orthogonal to $v$.

Let $K$ be a compact convex set in $\R^3$, and let $L$ be a plane containing the origin, and having $u$ as unit normal vector.
In this case there exists real concave functions $f$ and $g$ on $K|L$ such that $f(x)+g(x)\geq 0$ for $x\in K|L$, and
\begin{equation}
\label{bdKparameter}
K=\{x+tu:\, x\in K|L \mbox{ and }-g(x)\leq t\leq f(x)\}.
\end{equation}
Then the Steiner symmatrial $K_L$ of $K$ with respect to $L$ is
$$
K_L=\left\{x+tu:\, x\in K|L \mbox{ and }-\frac{f(x)+g(x)}2\leq t\leq \frac{f(x)+g(x)}2\right\}.
$$
Then $K_L$ is a convex convex set symmetric through $L$, which is a convex body of $K$ was a convex body. The following is a trivial observation that helps to understand the structure of the Steiner symmetrial of a polytope.

\begin{lemma}
\label{Steiner-linear}
Using the notation in (\ref{bdKparameter}),
if $\Pi\subset K|L$ is convex, and $f$ and $g$ are linear on $\Pi$, then $f+g$ is linear on $\Pi$.
\end{lemma}

For any compact $X$, we write ${\rm diam}\,X$ to denote the diameter of $X$.
Steiner symmetrial $K_L$ has the following basic properties (see Schneider \cite{Sch14}).
We write $B^3$ to denote the unit ball in $\R^3$ centered at the origin.

\begin{lemma}[Steiner]
\label{Steiner-symmetrization}
If $K$ is a convex body in $\R^3$ and $L$ is a linear $2$-plane, then
\begin{itemize}
\item $V(K_L)=V(K)$;
\item $S(K_L)\leq S(K)$ with equality if and only if $K$ is a translate of $K_L$;
\item $K\subset R\,B^3$ for $R>0$ implies $K_L\subset R\,B^3$.
\end{itemize}
\end{lemma}

The following observation ensures that the combinatorial type of a polytope is kept under Steiner symmetrization in certain situations.

\begin{lemma}
\label{possible-Steiner}
If the boundary of a polytope $P$ in $\R^3$ is triangulated using only the vertices of $P$ in a way such that two vertices $v_0$ and $v_1$ are connected by an edge to any other vertex, and any face contains either $v_0$ or $v_1$, 
and $L=(v_1-v_0)^\bot$, then $P_L$ is a double pyramid with the same numebr of vertices as $P$, and the induced triangulation of $P_L$  
  has the same combinatorial type as the one for $P$. 
\end{lemma}
\proof Let $\mathcal{T}$ be the family of triangles and edges in the triangulation of $P$.
If $T\in \mathcal{T}$ contains $v_i$ but not $v_{1-i}$ for $i\in\{0,1\}$, then let $s$ be the side of $T$ opposite to $v_i$. This edge $s$ of $T$ is the side of another triangle $T'\in \mathcal{T}$, and hence $T'$ contains $v_{1-i}$ but not $v_i$. We call $T'$ the twin of $T$, and we observe that $T$ is the unique twin of $T'$.

We observe that 
$$
\tilde{v}_1=(v_1|L)+\frac{v_1-v_0}2 \mbox{ \ and \ }
\tilde{v}_0=(v_1|L)+\frac{v_0-v_1}2
$$
 are vertices of $P_L$.
We define a bijection $\varphi$ between the faces of the triangulation of $P$ and the corresponding triangulation of $P_L$. If both $v_0$ and $v_1$ are vertices of a $T\in \mathcal{T}$, then $\varphi(T)=T_L$.

Next we assume that $T\in \mathcal{T}$ contains $v_i$ but not $v_{1-i}$ for $i\in\{0,1\}$, and $T$ is a proper subset of a face $F$ of $P$. We observe that 
 the common side $s$ of $T$ and $T'$ is a diagonal of $F$ as the endpoints of $s$ are vertices of $P$. Since
any vertex of $F$ is connected to either $v_i$ or $v_{1-i}$ by an edge, we deduce that $F$ is the quadrilateral
$T\cup T'$. Therefore $F_L$ is the quarilateral with diagonal $s|L$, and $s|L$ dissects $F_L$ into $\varphi(T)$ 
containing $\tilde{v}_i$ and $\varphi(T')$ 
containing $\tilde{v}_{1-i}$.

Finally we assume that $T\in \mathcal{T}$ contains $v_i$ but not $v_{1-i}$ for $i\in\{0,1\}$, and $T$ is a face 
of $P$, which we call a triangle of general type  of $\mathcal{T}$. In this case $T$ and its twin $T'$ are faces of a tetrahedron $S$. Any line $l$ parallel to the edge $[v_1,v_0]$
of $S$ intersecting ${\rm int}\,S$ meets ${\rm bd}\,S$ in one point of $T$ and in one point of $T'$, therefore
$$
\Pi=S|L=T|L=T'|L.
$$
We deduce from Lemma~\ref{Steiner-linear} that $P_L$ has two faces $\varphi(T)$ and $\varphi(T')$
where 
$$
\Pi=\varphi(T)|L=\varphi(T')|L,
$$
where $\varphi(T)$ 
contains $\tilde{v}_i$ and $\varphi(T')$ 
contains $\tilde{v}_{1-i}$.

Now any line $l$ parallel to $v_1-v_0$ intersecting ${\rm int}\,S$ and avoiding any edge in $\mathcal{T}$
meets ${\rm bd}\,P$ in one point of a triangle $T$ of general type and in one point of its twin $T'$. For the tetraeder $S$ determined by $T$ and $T'$, and for the edge $s=T\cap T'$, we have that $S_L$ is the tetrahedron 
determined by $s|L$, $\tilde{v}_1$ and $\tilde{v}_1$.
Therefore $\varphi$
is a combinatorial isomorphism of the triangulation of $P$ and the corresponding triangulation of $P_L$,
and each vertex of $P_L$ but $\tilde{v}_1$ and $\tilde{v}_1$ is contained in $L$.  
\proofbox

The basic properties of Steiner symmetrization yield directly the optimality of the regular tetrahedron among tetrahedra.
Before showing this, we point out a simple but useful observation that will be used several times in this paper.

\begin{lemma}
\label{inradius}
If a polytope $P$ in $\R^3$ contains a ball of radius $r$ touching each face of $P$, then
$$
\frac{S(P)^3}{V(P)^2}=\frac{27V(P)}{r^3}.
$$
\end{lemma}
\proof Writing $z$ to denote the center of the ball with radius $r$ in $P$, we dissect $P$ into pyramids whose apex 
is $z$ and base is a face of $P$. Calculating $V(P)$ as the sum of the volumes of these pyramids shows that
$V(P)=\frac13\,r S(P)$, which in turn yields Lemma~\ref{inradius}. \proofbox

\noindent{\bf Proof of Theorem~\ref{isoper-4vertices} } For any tetrahedron $T$, we write $r(T)$ to denote the radius inscribed ball touching each face of $T$.

For fixed $V>0$, let us consider a regular tetrahedron 
$P_0$ of volume $V$ containing the origin, and the family $\mathcal{P}$ of tetrahedrons $P$ of volume $V$ and containing the origin and satisfying
$S(P)\leq S(P_0)$. Let $z(P)$ be the center of the inscribed ball of $P\in \mathcal{T}$ where
Lemma~\ref{inradius} yields that $r(P)\geq r(P_0)$.

We claim that there exists $R>0$ depending on $V$ such that
\begin{equation}
\label{tetraRB}
P\subset RB^3\mbox{ \ for $P\in \mathcal{P}$.}
\end{equation}
 If $x\in P$, then either $\|z(P)\|\geq \|x\|/2$
or $\|x-z(P)||\geq \|x\|/2$ by the triangle inequality, therefore $P$ containes a right cone $C$ whose apex is either $z(P)$ or $o$, whose base is a circle of radius $r(P)$ and height at least $\|x\|/2$. We deduce that
$$
V(P)\geq V(C)\geq \frac{\|x\|/2}3\,r(P)^2\pi\geq  \frac{\|x\|}6\,r(P_0)^2\pi,
$$
completing the proof of (\ref{tetraRB}).

Since the volume and surface area are continuous function of the vertices of $P$, it follows from 
(\ref{tetraRB}) that there exists a $\widetilde{P}\in \mathcal{P}$ such that
$P=\widetilde{P}$ minimizes
$S(P)$. We suppose that $\widetilde{P}$ is not a regular tetrahedron, and seek a contradiction.
There exist vertices $v_1,v_2,v_3$ of $P$ such that
$\|v_3-v_1\|\neq \|v_3-v_2\|$. We use Steiner symmetrization with respect to
$L=(v_1-v_2)^\bot$. It follows from Lemma~\ref{Steiner-symmetrization} and
Lemma~\ref{possible-Steiner} that $\widetilde{P}_L\in\mathcal{T}$ and 
$S(\widetilde{P}_L)<S(\widetilde{P})$. This contradiction with the minimality of
$S(\widetilde{P})$ verifies Theorem~\ref{isoper-4vertices}.
\proofbox

To prove Proposition~\ref{likeoctahedron}, we write $|Q|$ to denote the area of a two-dimensional compact convex set $Q$.\\

\noindent{\bf Proof of Proposition~\ref{likeoctahedron} } Let $v_1,\ldots,v_6$ be the vertices of the polytope $P$ in a way such that each of them are of degree $4$, and $v_i$ and $v_{i+3}$ are not connected by an edge for $i=1,2,3$.
Let $L_1=(v_4-v_1)^\bot$, thus Lemma~\ref{possible-Steiner} applies
to $P'=P_{L_1}$. Let $v'_i$ be the vertex of $P'$ corresponding to $v_i$, $i=1,\ldots,6$.
In particular, $P'$ is a double pyramid with apices $v'_1$ and $v'_4$ over the the quadriteral with vertices $v'_2,v'_3,v'_5,v'_6\in L_1$, and $P'$ is symmetric through $L_1$. Next let 
$L_2=(v'_5-v'_2)^\bot$, thus Lemma~\ref{possible-Steiner} applies
to $P''=P'_{L_2}$. We write $v''_i$ to denote the vertex of $P''$ corresponding to $v_i$, $i=1,\ldots,6$.
As $L_1$ and $L_2$ are orthogonal, $P''$ is symmetric through both $L_1$ and $L_2$,
and $v''_3,v''_6\in L_1\cap L_2$. For $L_3=(v''_6-v''_3)^\bot$, the vertices of $P'''=P''_{L_3}$
are of the form $\pm t_iu_i$ for $i=1,2,3$ where $u_i$ is a unit normal to $L_i$ and $t_i>0$
for $i=1,2,3$. We deduce from Lemma~\ref{Steiner-symmetrization} that
\begin{equation}
\label{P'''}
V(P''')=V(P)\mbox{ \ and \ }S(P''')\leq S(P)
\end{equation}
with equality if and only if $P'''$ is a translate of $P$.

We observe $P'''$ can dissected into $8$ congruent copies of $[o,t_1u_1,t_2,u_2,t_3u_3]$,
and hence each face of $P'''$ is congruent to $F=[t_1u_1,t_2,u_2,t_3u_3]$ and
$$
V(P''')=\frac{4t_1t_2t_3}3\mbox{ \ and \ }S(P''')=8|F|.
$$
Let $r>0$ be the radius of largest ball centered at $o$ and contained in $P'''$, and let $u$ be the unit vector
satisfying that $ru$ is the closest point of $F$ to $o$. Since $\langle u,t_iu_i\rangle=r$ for $i=1,2,3$, we deduce
$$
ru=\sum_{i=1}^3\frac{r^2}{t_i}\,u_i=\sum_{i=1}^3\frac{r^2}{t_i^2}\,u_i.
$$
Here $ru\in [t_1u_1,t_2,u_2,t_3u_3]$ yields $\sum_{i=1}^3\frac{r^2}{t_i^2}=1$,
and we have $V(P''')=8\cdot\frac13\,r|F|=\frac13\,r S(P''')$. We conclude
using the Jensen inequality that
\begin{eqnarray*}
S(P)&\geq &S(P''')=\frac{3V(P''')}r=3V(P''')\sqrt{\sum_{i=1}^3\frac{1}{t_i^2}} 
\geq 3V(P''') 3^{\frac{1}2}(t_1t_2t_3)^{\frac{-1}3}=3^{\frac76}2^{\frac23}V(P''')^{\frac23}\\
&=& 
\sqrt[3]{108\sqrt{3}\,V(P)^2},
\end{eqnarray*}
with equality if and only if $P$ is a regular octahedron.
\proofbox

\section{If the combinatorial type is not of the octahedron}
\label{secnoocta}

We say that a polytope $P$ in $\R^3$ is strange with base $Q$ if $P$ is 
a double pyramid over a polygon $Q$ and $P$ has apices $v_1$ and $v_2$
where $v_1$ is a unit vector, $v_2=-v_1$, and
$$
Q=[o,v_3,v_4,v_5,v_6]\subset v_1^\bot;
$$
moreover either $Q$ is a triangle,  or
$o$ lies on the (relative) boundary of $Q$ and $v_3,v_4,v_5,v_6$ are vertices of $Q$ and $P$.
We note that
\begin{equation}
\label{volume}
V(P)=\mbox{$\frac23\,|Q|$}.
\end{equation}

The role of strange polytopes is explained by Lemma~\ref{non-octahedron} below. 
As a preparation for Lemma~\ref{non-octahedron},
a triangulated polytope $M$ in $\R^3$ means either a simplicial polytope, or a non-simplicial one whose boundary is triangulated using only the vertices of $M$. The boundary of $M$ is denoted by ${\rm bd}\, M$.
Now if $M$ is a triangulated polytope in $\R^3$ with vertices $w_1,\ldots,w_n$, and  $d_i$ is the degree of $w_i$ in the triangulation (the number of edges of the triangulation meeting
at $w_i$), then the Euler formula yields that
\begin{equation}
\label{degreesum}
\sum_{i=1}^n d_i=6n-12.
\end{equation}
We note that $3\leq d_i\leq n-1$ for $i=1,\ldots,n$.

\begin{lemma}
\label{non-octahedron}
If $P$ is a triangulated polytope with five or six vertices in $\R^3$ such that the triangulation is not combinatorially equivalent to the combinatorial type of the octahedron, then
some Steiner symmetrial of $P$ is similar to a strange polytope of the same number of vertices as $P$.
\end{lemma}
\proof If $P$ has five vertices  then (\ref{degreesum}) yields that any vertex is of degree either $3$ or $4$, 
two vertices, say $w_1$ and $w_2$ are of degree $3$, and three are of degree $4$; moreover, any face of the triangulation holds two vertices of degree (with respect to the triangulation) $4$, and one vertex of degree $3$. 
In particular, Lemma~\ref{possible-Steiner} applies for the Steiner symmetrization with
to the linear plane $L$ orthogonal to $w_1-w_2$, and hence $P_L$ is a translated and dilated image of
a strange polytope whose apices are $\pm\frac{w_1-w_2}{\|w_1-w_2\|}$.

Next, we assume that $P$ has six vertices. The combinatorial type of the octahedron is characterized by the facts that the polytope has six vertices, and each is of degree $4$. Thus if $P$ is a triangulated polytope with a six vertices in $\R^3$ such that the triangulation is not combinatorially equivalent to the octahedron's one, then it has a vertex $w_1$ of degree three by (\ref{degreesum}). Let $w_2,\ldots,w_6$ be the other vertices. Now $w_1$ is the end point of at least three edges of the polytope $P$. It follows that there exists three vertices, say $w_4,w_5,w_6$ connected to $w_1$ by an edge of $P$. Thus $[w_4,w_5,w_6]$ is a triangular face of the three-dimensional polytope 
$P_0=[w_2,\ldots,w_6]$. The consideration above show that for the triangulation  of ${\rm bd}P_0$ resulting from the triangulation of  ${\rm bd}P$, we may assume that $w_5$ and $w_6$ are degree $4$ vertices of 
the triangulation of ${\rm bd}P_0$.
Therefore $w_5$ and $w_6$ are degree $5$ vertices of the triangulation of  ${\rm bd}P$. 
In particular,  $w_5$ and $w_6$ are connected by an edge to any other vertex of $P$, and any face 
of $P$ contains either $w_5$ or $w_6$. We conclude from 
Lemma~\ref{possible-Steiner} that applying a Steiner symmetrization 
with respect to the linear plane $L$ orthogonal to $w_5-w_6$, $P_L$ is a translated and dilated image of
a strange polytope whose apices are $\pm\frac{w_5-w_6}{\|w_5-w_6\|}$, completing the proof of
Lemma~\ref{non-octahedron}.
\proofbox

We note that if $O$ is a regular octahedron, then
\begin{equation}
\label{188}
\frac{S(O)^3}{V(O)^2}=108\sqrt{3}<188,
\end{equation}
and if $D$ is the double pyramid of Theorem~\ref{isoper-5vertices}, then
\begin{equation}
\label{344}
\frac{S(D)^3}{V(D)^2}=243\sqrt{2}<344,
\end{equation}
The estimates of Lemmas~\ref{volumeest} and \ref{distanceest} show that a strange polytope with five or six vertices
should have a bounded shape.

\begin{lemma}
\label{volumeest}
Let $P$ be a strange polytope  with base $Q$.
\begin{description}
\item{(i)}
If either $V(P)\leq 0.274$ or
$V(P)\geq 3.4$; or equivalently, either $|Q|\leq 0.411$ or $|Q|\geq 5.1$, then
$$
\frac{S(P)^3}{V(P)^2}> 188.
$$
\item{(ii)}
If either $V(P)\leq 0.06$ or
$V(P)\geq 10$; or equivalently, either $|Q|\leq 0.09$ or $|Q|\geq 15$, then
$$
\frac{S(P)^3}{V(P)^2}> 344.
$$
\end{description}
\end{lemma}
\proof
Let $P_0$ be the Schwarz-rounding of $P$ with respect to the third coordinate axis. In particular, for any $t\in(-1,1)$, let
$H_t=\{(x,y,t):\,x,y\in\R\}$, then $H_t\cap P_0$ is a circular disk with center $(0,0,t)$ and area $|H_t\cap P|$.
It is known (see Schneider \cite{Sch14}) that $V(P_0)=V(P)$ and $S(P_0)\leq S(P)$. If $|Q|=r^2\pi$ for $r>0$, then
the boundary of $P_0$ can be unfolded into two sectors of the circular disk of radius $\sqrt{1+r^2}$, both bounded by an arc of length $2\pi r$. Therefore $S(P_0)=2r\pi \sqrt{1+r^2}$ and
$$
\frac{S(P_0)^3}{V(P_0)^2}=\frac{18\pi(1+r^2)^{\frac32}}{r}=18\pi\sqrt{\frac{(1+r^2)^3}{r^2}}
=18\sqrt{\frac{(\pi+|Q|)^3}{|Q|}}.
$$
Setting $s=|Q|$, we observe that the function
$$
f(s)=\frac{(\pi+s)^3}s=3\pi^2+\frac{\pi^3}s+3\pi s+s^2
$$
satisfies
$$
f'(s)=\frac{-\pi^3+3\pi s^2+2s^3}{s^2}=
\frac{(s+\pi)^2(2s-\pi)}{s^2}.
$$
In particular, $f(s)$ is strictly monotone decreasing for $0<s<\frac{\pi}2$, and strictly monotone increasing for $s>\frac{\pi}2$.
Since $18\sqrt{f(0.411)}>188$ and $18\sqrt{f(5.1)}>188$, 
and $18\sqrt{f(0.09)}>344$ and $18\sqrt{f(15)}>344$,
we conclude Lemma~\ref{volumeest}.
\proofbox

\begin{lemma}
\label{distanceest}
Let $P=[v_1,\ldots,v_6$ be a strange polytope where $v_1,v_2$ are the apices.
\begin{description}
\item{(i)}
If there exists some $i\in\{3,4,5,6\}$ such that $\|v_i\|\geq 6.5$, then
$$
\frac{S(P)^3}{V(P)^2}> 188.
$$
\item{(ii)}
If there exists some $i\in\{3,4,5,6\}$ such that $\|v_i\|\geq 17$, then
$$
\frac{S(P)^3}{V(P)^2}> 344.
$$
\end{description}
\end{lemma}
\proof Both for (i) and (ii), we observe that
$$
S(P)> S([v_i,v_5,v_6])= 2\|v_i\|.
$$
For (i),  we may assume that $V(P)<3.4$ by Lemma~\ref{volumeest}, and hence
$$
\frac{S(P)^3}{V(P)^2}> \frac{8\|v_i\|^3}{3.4^2}\geq  \frac{8\cdot 6.5^3}{3.4^2}>188.
$$
For (ii),  we may assume that $V(P)<10$ by Lemma~\ref{volumeest}, and hence
$$
\frac{S(P)^3}{V(P)^2}> \frac{8\|v_i\|^3}{10^2}\geq  \frac{8\cdot 17^3}{10^2}>344.
$$
\proofbox

\noindent{\bf Proof of Theorem~\ref{isoper-5vertices} } According to 
Theorem~\ref{isoper-4vertices}, we may assume that $P$ has five vertices.
It follows from Lemma~\ref{Steiner-symmetrization}, Lemma~\ref{non-octahedron} 
and (\ref{344}) that Theorem~\ref{isoper-5vertices} follows from the statement that
\begin{equation}
\label{strange5-claim}
\frac{S(P)^3}{V(P)^2}\geq 243\sqrt{2} \mbox{ \ if $P$ is a strange polytope with five vertices } 
\end{equation}
with equality if and only if  $P$ is congruent to the strange polytope $P_0$ 
whose the base $Q_0$ is a regular triangle whose insrcibed circle is centred at $o$ and has radius $\sqrt{2}$.
We note that $\frac{S(P_0)^3}{V(P_0)^2}= 243$ holds for the polytope $P_0$, the edges of $Q_0$ are of length $2\sqrt{6}$, and the other six edges of $P_0$ are of length $3$ (compare Theorem~\ref{isoper-5vertices}).

Let $\mathcal{P}$ be the family of strange polytopes $P$ such that the base $Q$ of $P$ is a triangle
and $\frac{S(P)^3}{V(P)^2}\leq 243$. It follows from Lemmas~\ref{volumeest} and \ref{distanceest}
that $|Q|\geq 0.18$ and $Q\subset 11B^3$.
Since the volume and surface area are continuous function of the vertices of $P$, we deduce the existence a $\widetilde{P}\in \mathcal{P}$ with base $\widetilde{Q}$ such that
$P=\widetilde{P}$ minimizes
$S(P)^3/V(P)^2$. We claim that 
\begin{equation}
\label{base5regular}
\mbox{$\widetilde{Q}$ is a regular triangle such that $o$ is the center of mass.}
\end{equation}
We suppose that (\ref{base5regular}) does not hold, and hence either $o$ is a vertex, or the vertices of $\widetilde{Q}$ are vertices of $\widetilde{P}$, as well, and we may denote them by $v_3,v_4,v_5$ in a way such that the perpendicular bisector of
$[v_3,v_4]$ misses $o$ or $v_5$. Now if $o$ is a vertex of $\widetilde{Q}$, then $\widetilde{P}$ has four vertices,
thus $\frac{S(\widetilde{P})^3}{V(\widetilde{P})^2}\geq 216\sqrt{3}$ by Theorem~\ref{isoper-4vertices},
contradicting $\frac{S(\widetilde{P})^3}{V(\widetilde{P})^2}\leq 243$.
Therefore the vertices $v_3,v_4,v_5$ of $\widetilde{Q}$ are vertices of $\widetilde{P}$ and  the perpendicular bisector of $[v_3,v_4]$ misses $o$ or $v_5$. We apply Steiner symmetrization with respect $L=(v_3-v_4)^\bot$,
thus Lemma~\ref{possible-Steiner} applies, and yields $\widetilde{P}_L\in \mathcal{P}$.
Since $v_1$ and $v_2$ are vertices of both $\widetilde{P}$ and $\widetilde{P}_L$,
we deduce from Lemma~\ref{Steiner-symmetrization} that
$S(\widetilde{P}_L)<S(\widetilde{P})$ and $V(\widetilde{P}_L)=V(\widetilde{P})$.
This contradicts the minimality of $\frac{S(\widetilde{P})^3}{V(\widetilde{P})^2}$, and proves 
(\ref{base5regular}). 

Let $\varrho$ be the inradius of $\widetilde{Q}$; namely, the common distance of $o$ from the sides of $\widetilde{Q}$.
It follows that the ball $rB^3$ is touched by each side of 
$\widetilde{P}$ for $r=\frac{\varrho}{\sqrt{1+\varrho^2}}$. We deduce from Lemma~\ref{inradius} 
and $|\widetilde{Q}|=3\sqrt{3}\,\varrho^2$ that
$$
\frac{S(\widetilde{P})^3}{V(\widetilde{P})^2}=\frac{27V(\widetilde{P})}{r^3}=
\frac{18|\widetilde{Q}|}{r^3}=54\sqrt{3}\,\frac{(1+\varrho^2)^{\frac32}}{\varrho}=
54\sqrt{3}\,\left(\frac{1+\varrho^2}{\varrho^{\frac23}}\right)^{\frac32}.
$$
It follows from differentiation that the unique place where  $\frac{1+\varrho^2}{\varrho^{\frac23}}$ attains its minimum is $\varrho=\sqrt{2}$, yielding   (\ref{strange5-claim}).
\proofbox

\noindent{\bf Proof of Theorem~\ref{isoper-6vertices} } According to 
Theorem~\ref{isoper-5vertices}, we may assume that the polytope $P$ has six vertices.

It follows from Proposition~\ref{likeoctahedron}, Lemma~\ref{Steiner-symmetrization}, Lemma~\ref{non-octahedron} and (\ref{188}) that Theorem~\ref{isoper-6vertices} follows from the statement that
\begin{equation}
\label{strange6-claim}
\frac{S(P)^3}{V(P)^2}\geq 188 \mbox{ \ if $P$ is a strange polytope. } 
\end{equation}
According to Lemma~\ref{volumeest} and Lemma~\ref{distanceest}, we may assume that the strange polytope $P$ with base $Q$ satisfies
\begin{description}
\item{(a)} $0.411 \leq |Q|\leq 5.1$
\item{(b)} $P\subset 6.5 B^3$
\end{description}
As described at the beginning of this section, the strange polytope $P$ is a double pyramid with apices $v_1$ and $v_2$
where $v_1$ is a unit vector, $v_2=-v_1$, and
$$
Q=[o,v_3,v_4,v_5,v_6]\subset v_1^\bot
$$
where $o$ lies on the (relative) boundary of $Q$ and $v_3,v_4,v_5,v_6$ are vertices of $Q$ and $P$. We may also assume that
$o,v_3,v_4,v_5,v_6$ lie in this order on the  (relative) boundary of $Q$.
In particular $P$  depends on $7$ real parameters; namely,  $v_3,\ldots,v_6$ have $8$ coordinates in $v_1^\bot$, but we can fix ${\rm lin}\,v_6$. 

In order to simplify calculations; more precisely, to reduce the number of free parameters, we apply one more Steiner symmetrization to $P$ with respect to $L=(v_3-v_6)^\bot$. Set $l={\rm lin}(v_3-v_6)$. In order to descrive $Q_L$,
possible after reversing the order of $v_3,\ldots,v_6$, we may assume that $v_4$ is not further from $l$, than $v_5$.
It follows from Lemma~\ref{Steiner-linear} that 
$$
Q_L=[o,w_1,w_2,w_3,w_4,w_5]
$$
where $w_1=\frac12(v_3-v_6)$ and $w_5=\frac{-1}2(v_3-v_6)$, $w_2-w_4=\lambda(v_3-v_6)$ for $\lambda>0$
and $w_3\in L\cap{\rm lin}Q$. Here $w_2,w_4\in v_4+l$ and $w_3\in v_5+l$, $o$ is a vertex of $Q_L$ if and only if it is a vertex of $Q$, and $w_3\in[w_2,w_4]$ if and only if $v_3-v_6$ is parallel to $v_4-v_5$.
We choose an orthonomal basis of $\R^3$ where $w_3/\|w_3\|$ is the first basis vector, 
$(w_2-w_5))/\|w_2-w_5)\|$ is the second basis vector, and $v_1$ is the third basis vector.

With respect to this basis, we have $v_1=(0,0,1)$ and $v_2=(0,0,-1)$, and for suitable
\begin{equation}
\label{cond1}
0\leq x_1\leq x_2\leq x_3\mbox{ \ and \ }y_1,y_2\geq 0,
\end{equation}
we have $w_1=(x_1,y_1,0)$, $w_2=(x_2,y_2,0)$, $w_3=(x_3,0,0)$, $w_4=(x_2,-y_2,0)$
and $w_5=(x_1,-y_1,0)$. Since $w_1,w_2,w_3$ lie in this order on the relative boundary of $Q_L$, we deduce that
\begin{equation}
\label{cond2}
x_2y_1-x_1y_2\geq 0\mbox{ \ and \ }
(x_3-x_1)y_2-(x_3-x_2)y_1\geq 0.
\end{equation}
We observe that $|Q_L|=|Q|$ where
$$
|Q_L|=x_2y_1-x_1y_2+x_3y_2,
$$
and hence (a) yields that
\begin{equation}
\label{cond3}
x_2y_1-x_1y_2+x_3y_2\geq 0.411.
\end{equation}
In addition, (b) implies that
\begin{equation}
\label{cond4}
x_3,y_1,y_2\leq 6.5.
\end{equation}
We note that $V(P_L)=\frac23\,|Q_L|$ and
$$
S(P_L)=2\|v_1\times w_1\|+2\|(w_1-v_1)\times(w_2-v_1)\|+2\|(w_2-v_1)\times(w_3-v_1)\|.
$$
As functions of $x_1,x_2,x_3,y_1,y_2$, we consider the functions $S$ and $V$ defined by
\begin{eqnarray}
\label{sdef}
S=S(P_L)&:=&
2\sqrt{x_1^2+y_1^2}+2\sqrt{ (x_1-x_2)^2+(y_1-y_2)^2+(x_2y_1-x_1y_2)^2 }+\\
\nonumber
&& +2\sqrt{ y_2^2+(x_3-x_2)^2+x_3^2y_2^2};\\
\label{vdef}
V=V(P_L)&:=&\frac23\,(x_2y_1-x_1y_2+x_3y_2).
\end{eqnarray}

We have checked using Wolfram Mathematica that
\begin{equation}
\label{mutant6-claim}
S^3-188 \cdot V^2>3.44
\end{equation}
holds for the functions $S$ and $V$ defined in (\ref{sdef}) and (\ref{vdef}) under the conditions
(\ref{cond1}), (\ref{cond2}), (\ref{cond3}) and (\ref{cond4}) for the variables $x_1,x_2,x_3,y_1,y_2$.
Therefore, we conclude Theorem~\ref{isoper-6vertices} from
 Lemma~\ref{Steiner-symmetrization}, (\ref{strange6-claim}) and (\ref{mutant6-claim}).  \proofbox

\end{document}